\documentclass[12pt]{article}
\usepackage{amsmath,amsfonts,amssymb,amscd}
\usepackage{graphicx}

\newtheorem{theo}{Theorem}

\newtheorem{prop}{Proposition}[section]

\newtheorem{cor}{Corollary}[section]
\newtheorem{rem}{Remark}[section]

\makeatletter \@addtoreset{equation}{section} \makeatother

\newcommand{\mR}{\mathbb{R}}

\newcommand{\mZ}{\mathbb{Z}}
\newcommand{\mN}{\mathbb{N}}

\newcommand{\bW}{{\bf W}}

\newcommand{\bb}{{\bf b}}

\newcommand{\be}{{\bf e}}

\newcommand{\bh}{{\bf h}}

\newcommand{\calA}{{\cal A}}

\newcommand{\calI}{{\cal I}}
\newcommand{\calJ}{{\cal J}}
\newcommand{\calK}{{\cal K}}

\newcommand{\calP}{{\cal P}}

\newcommand{\calT}{{\cal T}}

\newcommand{\calX}{{\cal X}}

\newcommand{\eps}{\varepsilon}
\newcommand{\ph}{\varphi}

\newcommand{\Deg}{\operatorname{Deg}}

\newcommand{\one}{{\bf{1}}}
\newcommand{\Par}{\operatorname{Par}}

\newcommand{\symm}{\operatorname{symm}}

\begin{document}

%\large
\title
{Isochronicity in 1 DOF}
\author{ Dmitry Treschev \\
Steklov Mathematical Institute of Russian Academy of Sciences
%\footnote{This work was performed at the Steklov International Mathematical Center and supported by the %Ministry of Science and Higher Education of the Russian Federation (agreement no. 075-15-2019-1614).}
}
\date{}
\maketitle

\begin{abstract}
Our main result is the complete set of explicit conditions necessary and sufficient for isochronicity of a Hamiltonian system with one degree of freedom. The conditions are presented in terms of Taylor coefficients of the Hamiltonian function.
\end{abstract}

\section{Introduction}

\subsection{Main question}

Consider the Hamiltonian system
\begin{equation}
\label{ham_sys}
  \dot x = \partial H / \partial y, \quad
  \dot y = - \partial H / \partial x,
\end{equation}
where $(x,y)$ are coordinates on the plane $\mR^2$ and $H=H(x,y)$ is a smooth function. We assume that
\begin{equation}
\label{HHH}
  H = H_2 + H_*, \qquad
  H_2 = \pi (x^2 + y^2), \quad
  H_* = O(|x|+|y|)^3.
\end{equation}
Then the point $(x,y)=(0,0)$ is an equilibrium position and all solutions of the linearized system (the system with Hamiltonian $H_2$) are periodic with period 1. Our main question is as follows.
\smallskip

{\it For which functions $H_*$ all solutions are periodic with period one?}
\smallskip

Such systems are said to be {\it isochronous}.

\subsection{Main theorem}

It is convenient to use the complex coordinates $(z,\overline z)$
$$
  z = \frac1{\sqrt2} (x + iy), \quad
  \overline z = \frac1{\sqrt2} (x - iy).
$$
Then
\begin{equation}
\label{H2=zz}
  H_2 = 2\pi z\overline z.
\end{equation}

Let
\begin{equation}
\label{H=sumH}
  H_* = \sum_{|\beta|\ge 3} H_\beta \mu^\beta , \qquad
  \mu^\beta = z^{\beta_1} \overline z^{\beta_2}, \quad
  |\beta| = \beta_1 + \beta_2, \quad
  \beta\in\mZ_+^2,
\end{equation}
where $\mZ_+ = \{0,1,2,\ldots\}$. We put
\begin{equation}
\label{Sigmadk}
    \Sigma^d_k
  = \frac1{(2\pi)^k} \sum  H_{\beta^{(1)}}\cdot\ldots\cdot H_{\beta^{(k)}},
\end{equation}
where the summation is performed under the conditions
\begin{equation}
\label{cond_Sigmadk}
  \beta^{(1)},\ldots,\beta^{(k)}\in\mZ_+^2, \quad
  |\beta^{(j)}| \ge 3, \quad
  \beta^{(1)}+\cdots+\beta^{(k)} = (d+k,d+k) .
\end{equation}

\begin{rem}
\label{rem:1kd}
Conditions (\ref{cond_Sigmadk}) imply
$$
  3k\le |\beta^{(1)} + \ldots + \beta^{(k)}| = |(d+k,d+k)| = 2(d+k).
$$
Therefore, if the sum $\Sigma_d^k$ contains at least one term, then $1\le k\le 2d$.
\end{rem}

\begin{theo}
\label{theo:main}
Suppose $H_2$ satisfies (\ref{H2=zz}), and the series (\ref{H=sumH}) converges in a neighborhood of the point $(z,\overline z)=(0,0)$. Then the system (\ref{ham_sys}) is (locally) isochronous if and only if the following infinite collection of equations hold:
\begin{equation}
\label{eq_conj}
  \sum_{s=1}^d \frac{(-1)^{s-1} (d+s)!}{s! (d+1)!} \Sigma_s^d = 0, \qquad
  d = 1,2,\ldots
\end{equation}
\end{theo}

\subsection{Comments}

1. The problem of isochronicity (or superintegrability) has a long history going back to Huygens who constructed an isochronous cycloidal pendulum. More general results on the existence of isochronous systems are established for some special classes of Hamiltonian systems (\ref{ham_sys}). One of such classes is the ``natural'' or ``physical'' one, where $H=y^2/2+V(x)$. Here the function $V$ is said to be a potential. Here we mention results, presented in \cite{Z13}, concerning construction of isochronous potentials, which generalize all known examples by Piskunov, Urabe, Stillinger, Dorignac, etc. Note that unlike our assumptions the function $V$ in \cite{Z13} may be only $C^1$-smooth.

In  \cite{Still,BM2003} some results on the existence of globally isochronous systems with 1 degree of freedom are presented. The monograph \cite{Calogero_book} contains several constructions which work in the multidimensional case as well.
\medskip

2. Isochronicity condition means that the normal form of $H$ equals $H_2$. Equivalently for some generating function $z\overline w + S(z,\overline w)$, we have:
$$
  H(z,\overline w + \partial_z S) = H_2(z + \partial_{\overline w} S, \overline w).
$$
At first glance study of this equation looks simpler than the method we choose below. But we do not know how our main result (equations (\ref{eq_conj})) can be obtained in this way.
\medskip

3. Equations (\ref{eq_conj}) may be used to solve the following class of problems. Suppose that the Hamiltonian (\ref{HHH}) contains as a parameter a function $\ph$ of one variable:
$$
  H = H^\ph(x,y) = \pi(x^2 + y^2) + O(|x|+|y|)^3.
$$
If the family of functions $H^\ph$ is chosen reasonably, Taylor coefficients of $\ph$ are computed uniquely from equations (\ref{eq_conj}).

If we plug such $\ph$ in $H^\ph$, we obtain an isochronous system provided the Taylor series, that determines $H^\ph$, converges. Convergence or divergence of this series is a separate problem.
\medskip

{\bf Example 1}. Suppose that
$H =  2\pi z\overline z
    + az^4 + bz^3\overline z + \overline bz\overline z^3 + \overline a\overline z^4
    + \ph(z\overline z)$,
\begin{equation}
\label{ph}
     \ph(w)
  = \sum_{j=3}^\infty \ph_j w^j.
\end{equation}
In this case conditions (\ref{eq_conj}) are satisfied for a unique series (\ref{ph}). Numeric study gives an evidence that this series diverges.
\medskip

{\bf Example 2}. N. Elfimov, \cite{Elf}
$$
  H = \ph(y^2) f(x^2), \qquad
  \ph(0) = f(0) = 1, \quad
  \ph'(0) = f'(0) = 1.
$$

\begin{prop}
For any real-analytic $f$ there exists a real-analytic $\ph$ such that the system (\ref{ham_sys}) is isochronous.
\end{prop}

{\bf Example 3}. An analogous convergence problem combined with superintegrability phenomenon appears in a slightly another context. Consider a billiard system on a billiard table with a smooth boundary, having two perpendicular symmetry axes. Such a system has two trajectories with period 2. Take one of them.

In \cite{Tre_PhD} we ask the following question. {\it Do the above tables exist such that locally near this trajectory the square of the billiard map is conjugated to a rigid rotation on a plane?}

Results (partially analytic, partially numeric), essentially confirming the positive answer to the question are contained in \cite{Tre_PhD,Tre_2015,Tre_2018}. A multidimensional version of this problem is discussed in \cite{Tre_2017}.

In this problem the boundary curve plays the role of the function $\ph$. Under reasonable nonresonance conditions superintegrability (linearizability) of the billiard map fixes uniquely the Taylor expansion of $\ph$ (up to a real parameter). So, the main problem is the convergence of the corresponding series. This problem remains unsolved.

\section{Sketch of the proof of Theorem \ref{theo:main}}

\subsection{The main equation}

We introduce the Poisson bracket
$$
    \{f,g\}
  = i \Big(\frac{\partial f}{\partial\overline z_j} \frac{\partial g}{\partial z_j}
          - \frac{\partial f}{\partial z_j} \frac{\partial g}{\partial\overline z_j}\Big) \quad
    \mbox{for any two functions $f=f(z,\overline z)$ and $g=g(z,\overline z)$}.
$$

For any function $F = F(z,\overline z)$ let $\widehat F$ denote the operator
$$
  \ph\mapsto\widehat F\ph := \{F,\ph\}, \qquad
  \ph = \ph(z,\overline z).
$$
According to the Jacobi identity for any two functions $F',F''$ the commutator
$[\widehat F',\widehat F'']$ satisfies the equation
$$
  [\widehat F',\widehat F''] = \widehat F, \qquad
  F = \{F',F''\}.
$$

For any function $\ph : U\to\mR$ we have $\partial_t\ph = \widehat H\ph$. Hence, the flow $g^t$ of the system (\ref{ham_sys}) has the form
$$
  g^t = e^{t\widehat H}
      = 1 + \frac{t}{1!}\widehat H + \frac{t^2}{2!}\widehat H^2 + \ldots ,
$$
where $1$ denotes the identity operator. Isochronicity of the system (\ref{ham_sys}) takes the form $e^{\widehat H} = 1$. Below we use the more detailed form of this equation:
$$
    \bW = 0, \qquad
    \bW = 1 + \frac1{1!} \big(\widehat H_2 + \widehat H_*\big)
            + \frac1{2!} \big(\widehat H_2 + \widehat H_*\big)^2
            + \ldots
            - 1 .
$$

\subsection{Step 1: gather terms in $\bW$}

The operator $\bW$ is a linear combination of operators
$\widehat F_1 \widehat F_2 \ldots \widehat F_k$, where any operator $\widehat F_j$ equals either $\widehat H_2$ or $\widehat H_*$. Below we use the expansion (\ref{H=sumH}). Therefore we will deal with (non-commutative) products of the following form: $\chi=\widehat F_1 \widehat F_2 \ldots \widehat F_k$, where any operator $\widehat F_j$ equals either $\widehat H_2$ or $\widehat\mu^\beta$, $|\beta|\ge 3$. We define
$$
    (\mZ_+^2)^k_*
  = \{\bb = (\beta^{(1)},\ldots,\beta^{(k)})
               : \beta^{(j)}\in\mZ_+^2, \, |\beta^{(j)}|\ge 3\}.
$$

Let $\calX_{\bb,s}$, $\bb\in (\mZ_+^2)^k_*$, $s\ge k$ be the set of all products $\chi$ such that for some sequence $1\le j_1<j_2<\ldots < j_k\le s$ we have:
$\widehat F_{j_1}=\widehat\mu^{\beta^{(1)}}, \ldots,\widehat F_{j_k}=\widehat\mu^{\beta^{(k)}}$ and all other multipliers $\widehat F_j$ equal $\widehat H_2$.

In the other words, any product $\chi\in\calX_{\bb,s}$ is obtained from $\widehat\mu^{\beta^{(1)}}\ldots\widehat\mu^{\beta^{(k)}}$ by inserting in arbitrary places $s-k$ copies of the operator $\widehat H_2$.

The operator $\bW$ takes the form
\begin{equation}
\label{W_no1}
    \bW
  = \sum_{0\le k\le s<\infty} \frac1{s!}
     \sum_{\bb\in(\mZ_+^2)_*^k}
      \sum_{\chi\in\calX_{\bb,s}} H^\bb \chi
    - 1, \qquad
    H^\bb
  = H_{\beta^{(1)}} H_{\beta^{(2)}} \ldots H_{\beta^{(k)}}.
\end{equation}

For any $\bb\in (\mZ_+^2)_*^k$ we put
\begin{equation}
\label{S_m}
    \widehat S_\bb
  = \sum_{s=k}^\infty \frac1{s!} \sum_{\chi\in\calX_{\bb,s}} \chi .
\end{equation}
Then by (\ref{W_no1})
\begin{equation}
\label{W_no2}
    \bW
  = \sum_{k=0}^\infty \sum_{\bb\in (\mZ_+^2)_*^k}  H^\bb\widehat S_\bb
    - 1.
\end{equation}

\subsection{Step 2: computation of $\widehat S_\bb$}

For any function $F=F(z,\overline z)$ we use the notation $\bh\widehat F = [\widehat H_2,\widehat F]$.
Then $\bh\widehat F = \widehat G$, $G=\{H_2,F\}$.

Our further idea is to move in any monomial $\chi$ the multipliers $\widehat H_2$ to the right by using the identities
$\widehat H_2 \bh^s\widehat H_n = \bh^s\widehat H_n \widehat H_2 + \bh^{s+1}\widehat H_n$.
To present a formula for $\widehat S_\bb$, we need another notation.

For any holomorphic at the origin function $\Phi=\Phi(\xi_1,\ldots,\xi_k)$,
$$
    \Phi(\xi)
  = \sum_{\alpha_1,\ldots,\alpha_k\ge 0}
             \Phi_{\alpha_1,\ldots,\alpha_k} \xi_1^{\alpha_1}\ldots\xi_k^{\alpha_k} ,
$$
we define the $k$-linear operator-valued form $\widehat\Phi$, acting on holomorphic at the origin functions $F_1,\ldots,F_k$:
$$
     \widehat\Phi(F_1,\ldots,F_k)
  := \sum_{\alpha_1,\ldots,\alpha_k\ge 0}
        \Phi_{\alpha_1,\ldots,\alpha_k} \bh^{\alpha_1}\widehat F_1\ldots\widehat\bh^{\alpha_k} F_k.
$$

For any two natural numbers $m\le n$ and $\eps\in\mR$ we define $\rho^-_{m,{m-1}}=\rho^+_{m,{m-1}}=1$ and
\begin{eqnarray*}
     \rho^+_{m,n,\eps}
 &=& (\eps + \xi_m)(\eps + \xi_m + \xi_{m+1})
             \ldots (\eps + \xi_m + \xi_{m+1} + \ldots + \xi_n),  \\
     \rho^-_{m,n,\eps}
 &=& (-1)^{n-m+1} (\xi_m + \ldots + \xi_{n-1} + \xi_n + \eps)
          \ldots  (\xi_{n-1} + \xi_n + \eps)(\xi_n + \eps) , \\
\end{eqnarray*}

\begin{prop}
\label{prop:S_m}
$\widehat S_\emptyset = 1$ and for any
$\bb = (\beta^{(1)},\ldots,\beta^{(k)})\in (\mZ_+^2)_*^k$, $k\ge 1$
\begin{equation}
\label{SandPhi}
     \widehat S_\bb
  =  \widehat\Phi_s (\mu^{\beta^{(1)}},\ldots,\mu^{\beta^{(k)}}), \qquad
          \Phi_k(\xi)
 :=       \sum_{\alpha_1,\ldots,\alpha_k\ge 0}
             \frac{\xi_1^{\alpha_1}\ldots \xi_k^{\alpha_k}}
                  {\calA_k(\alpha)} ,
\end{equation}
where
\begin{equation}
\label{Phi=sum}
    \Phi_k(\xi)
  =   \frac{e^{\xi_1+\ldots+\xi_k} - 1}{\rho^-_{1,0,0} \rho^+_{1,k,0}}
    + \frac{e^{\xi_2+\ldots+\xi_k} - 1}{\rho^-_{1,1,0} \rho^+_{2,k,0}}
    + \ldots
    + \frac{e^{\xi_k} - 1}{\rho^-_{1,k-1,0} \rho^+_{k,k,0}} .
\end{equation}
\end{prop}

\begin{cor}
\begin{equation}
\label{W_no3}
  \bW = \sum_{k=1}^\infty \sum_{\bb\in (\mZ_+^2)_*^k} H^\bb\widehat\Phi_k(\mu_\bb), \qquad
  \mu_\bb = (\mu^{\beta^{(1)}},\ldots,\mu^{\beta^{(k)}}).
\end{equation}
\end{cor}

\subsection{Step 3: computation of $\Phi_k$ at resonances}

The equations
$$
  \widehat H_2 \mu^\beta = \lambda_\beta \mu^\beta, \qquad
  \lambda_\beta = 2\pi i (\beta_1 - \beta_2)
$$
imply $\bh^s\widehat\mu^\beta = \lambda_\beta^s \widehat\mu^\beta$. Therefore
\begin{equation}
\label{Phi=Phi}
    \widehat\Phi_k(\mu^{\beta^{(1)}},\ldots,\mu^{\beta^{(k)}})
  = \Phi_k(\lambda_{\beta^{(1)}},\ldots,\lambda_{\beta^{(k)}})\,
       \widehat\mu^{\beta^{(1)}}\ldots\widehat\mu^{\beta^{(k)}}.
\end{equation}
Hence, we have to compute values of the functions $\Phi_k$ at points
$\xi=(\lambda_{\beta^{(1)}},\ldots,\lambda_{\beta^{(k)}})$ which lie in the set $2\pi i\mZ^k$.

Numerators of the fractions (\ref{Phi=sum}) vanish on the set $2\pi i\mZ^k$. Therefore $\Phi_k(\xi)=0$ if the numbers $\rho^-_{1,l-1,0} \rho^+_{l,k,0}$ do not vanish. The point $\xi\in 2\pi i\mZ^k$ at which one or several polynomials $\rho^\pm_{m,n}$ vanish, is naturally associated with a resonance.

For any set $K\subset\{1,\ldots,k\}$ consider the linear form
\begin{equation}
\label{ph_K}
    \ph_K(\xi)
  = \sum_{j\in K} \xi_j, \qquad
    \xi\in 2\pi i\mZ^k.
\end{equation}
Such forms will be of a special importance for us because the homogeneous polynomials
$\rho^-_{1,l-1,0} \rho^+_{l,k,0}$ are products of such forms.

For any point $\xi\in 2\pi i\mZ$ we have the set of resonances
$$
    \calK_\xi
  = \{K\subset\{1,\ldots,k\} : \ph_K(\xi)=0\}.
$$
We say that $\calK_\xi$ is generated by the sets $K_1,\ldots,K_m$ if
\begin{itemize}
\item the forms $\ph_{K_1},\ldots,\ph_{K_m}$, are linear independent and
\item for any $K\in\calK_\xi$ the forms $\ph_K,\ph_{K_1},\ldots,\ph_{K_m}$ are linear dependent.
\end{itemize}

We say that $K\subset \{1,\ldots,k\}$ is an interval if it has the form
\begin{equation}
\label{int}
  K = I(k^-,k^+) := \{j\in\mZ: k^-\le j\le k^+\},\qquad
  1\le k^-\le k^+\le k.
\end{equation}

First, consider the ``simple'' case when $\calK_\xi$ is generated by the intervals $I_1,\ldots,I_m$,
\begin{equation}
\label{Ks}
  I_s = I(k_s^-,k_s^+), \quad
  1\le k_s^-\le k_s^+ \le k, \qquad
  s = 1,\ldots,m ,
\end{equation}
where $k_{s+1}^- = k_s^+ + 1$ for any $s = 1,\ldots,m-1$.

This condition means that the intervals $I_1,\ldots,I_m$ pairwise do not intersect and their union is an interval. Such a collection of intervals will be said to be {\it friendly}.

\begin{prop}
\label{prop:Phi=}
Suppose that $\calK_\xi$ is generated by a friendly collection of intervals (\ref{Ks}). Then
\begin{eqnarray}
\nonumber
&\displaystyle
     \Phi_k(\xi)
  =  \frac1{(m-1)!}\partial_x^{m-1} \big|_{x=0} \be(x) N(x), & \\
\label{N(x)}
&\displaystyle
     \be(x)
  =  \frac{e^x - 1}{x}, \quad
     N(x)
  = \frac1{\rho^-_{1,k_1^- - 1,-x}
            \rho^-_{k_1^- + 1,k_1^+,-x}
             \ldots
              \rho^-_{k_m^- + 1,k_m^+,-x}
               \rho^+_{k_m^+ + 1,k,x}} . &
\end{eqnarray}
\end{prop}

The case of a general resonance may be reduced to a combination of ``simple'' ones.
Given a point $\xi\in 2\pi i\mZ^k$ let $\calI$ be the set of all intervals $I\subset\{1,\ldots,k\}$ such that $\ph_I(\xi)=0$.
We say that the subset $\calJ\subset\calI$ is minimal if any $K\in\calI\setminus\calJ$ is a union of several (more than one) intervals from $\calJ$ and any $K\in\calJ$ is not a union of several intervals from $\calJ$.
Obviously for any $\xi\in 2\pi i\mZ^k$ the minimal set of intervals exists and is unique.

The friendly collection $\calJ' = \{K_{j_1},\ldots,K_{j_s}\}\subset\calJ$ is said to be complete if there is no larger friendly collection $\calJ''\supset\calJ'$.
Complete friendly collections $\calJ^{(1)},\ldots,\calJ^{(\tau)}$ break the minimal set $\calJ$ into nonintersecting subsets:
$$
  \calJ = \cup_{j=1}^\tau \calJ^{(j)}, \qquad
  \calJ^{(j')}\cap\calJ^{(j'')} = \emptyset
  \mbox{ for } j'\ne j''.
$$

For any $j=1,\ldots,\tau$ we define the function $N_j$ by (\ref{N(x)}), where $m=m_j$ and the numbers $k_1^\pm,\ldots,k_{m_j}^\pm$ depend on $j$.

\begin{prop}
\label{prop:Phi_k}
For any $\xi\in 2\pi i\mZ^k$
\begin{equation}
\label{Phi=sum+}
    \Phi_k(\xi)
  = \sum_{j=1}^\tau \frac1{(m_j-1)!}
               \partial_x^{m-1}\big|_{x=0} \be(x) N_j(x) .
\end{equation}
\end{prop}

\subsection{Step 4: symmetrization of $\widehat\Phi_k$}

The group $S_k$ acts on $(\mZ_+^2)^k_*$ by permutations. For any $\sigma\in S_k$ we have
$$
          \bb = (\beta^{(1)},\ldots,\beta^{(k)})
  \mapsto \sigma(\bb) = (\beta^{(\sigma(1))},\ldots,\beta^{(\sigma(k))}).
$$

For any $\bb\in (\mZ_+^2)^k$ we introduce the notation
$$
    H^\bb
  = H_{\beta^{(1)}} \ldots H_{\beta^{(k)}} , \quad
    \mu_\bb
  = (\mu^{\beta^{(1)}},\ldots,\mu^{\beta^{(k)}}) , \quad
    \Psi_k(\mu_\bb)
  = \sum_{\sigma\in S_k} \Phi_k(\mu_{\sigma(\bb)}) .
$$
Then by (\ref{Phi=Phi})
$$
     \widehat\Phi_k(\mu_\bb)
  =  \sum_{\sigma\in S_k}\sum_{\bb\in (\mZ_+^2)^k_*}
            \frac{1}{k!} H^\bb \widehat\Phi_k(\mu_{\sigma(\bb)})
  =  \sum_{\bb\in (\mZ_+^2)^k_*}
            \frac{1}{k!} H^\bb \widehat\Psi_k(\mu_\bb).
$$

Hence
$$
    \bW
 =  \sum_{k=1}^\infty
     \sum_{\bb\in (\mZ_+^2)^k_*}
      \frac{1}{k!} H^\bb\widehat\Psi_k(\mu_\bb).
$$

\subsection{Step 5: computation of $\widehat\Psi_k$ at a resonance}

We use the notation $\bb = (\beta^{(1)},\ldots,\beta^{(k)})\in(\mZ_+^2)^k_*$
\begin{eqnarray*}
& \one = (1,1)\in\mZ_+^2, \quad
  \beta'\wedge\beta'' = \beta'_1\beta''_2 - \beta'_2\beta''_1, \qquad
  \beta',\beta''\in\mZ^2_+, & \\
& \langle\bb\rangle = \sum_{j=1}^k (\beta^{(j)} - \one), \quad
  \langle\bb\rangle_K = \sum_{j\in K} (\beta^{(j)} - \one), \qquad
  K\subset\{1,\ldots,k\}. &
\end{eqnarray*}

Possible resonances for $\widehat\Psi_k$ are $\ph_K(\xi)=0$, see (\ref{ph_K}).
These resonance conditions are equivalent to
$$
  \langle\bb\rangle_K = (d + \# K)\cdot\one \quad
  \mbox{for some}\quad
  d = d_{K,\bb}\in\mN.
$$

We define the set $\Par_n^{k,\kappa}$ of prepartitions of $\{1,\ldots,k\}$: we say that $\calT\in\Par_n^{k,\kappa}$ iff
\begin{itemize}
\item $\calT = \{A'_1,A''_1,\ldots,A'_n,A''_n\}$, $A'_s,A''_s\subset\{1,\ldots,k\}$,
\item intersection of any two sets from $\calT$ is empty,
\item $\# \cup_{s=1}^n (A'_s\cup A''_s) = \kappa$.
\end{itemize}

In particular, if $\kappa=k$ then $\cup_{s=1}^n (A'_s\cup A''_s) = \{1,\ldots,k\}$ and $\calT$ is a partition. We say that $\bb\in (\mZ_+^2)^k_*$ is compatible with $\calT\in\Par_n^{k,\kappa}$ (the notation is $\bb\vdash\calT$ or $\calT\dashv\bb$) if
$$
  \one\wedge\langle\bb\rangle_{A'_s} = 0 \quad
  \mbox{for any }  s = 1,\ldots,n.
$$

For any $\bb\in(\mZ_+^2)^k_*$ and $\{K,L\}\in\Par_1^{k,\kappa}$ such that $\bb\vdash\{K,L\}$ we put
\begin{equation}
\label{PKL}
    \calP_{K,L,\bb}
  = \alpha_{K,L,\bb}\, \widehat\mu^{\langle\bb\rangle_{K\cup L} + \one} , \quad
    \alpha_{K,L,\bb}
  = \frac{(-1)^{\# K - 1} d_{K,\bb}!}
         {(d_{K,\bb} - \# K - \# L + 1)! (2\pi)^{\# K + \# L - 1}} ,
\end{equation}
where $\langle\bb\rangle_K = (d_{K,\bb} + \# K)\cdot\one$.
\smallskip

For any collection of operators $\widehat F_1,\ldots,\widehat F_n$ we define
$$
     \symm(\widehat F_1,\ldots,\widehat F_n)
  =  \frac1{n!} \sum_{\sigma\in S_n}
       \widehat F_{\sigma(1)} \cdots \widehat F_{\sigma(n)} .
$$
We associate with any pair $\bb\vdash\calT\in\Par_n^{k,\kappa}$ the operator
\begin{equation}
\label{ABAB}
    \calP_{\calT,\bb}
  = \symm \big( \calP_{A'_1,A''_1,\bb}, \ldots , \calP_{A'_n,A''_n,\bb} \big) .
\end{equation}

\begin{prop}
\label{prop:Psi(mu)}
For any $\bb\in (\mZ_+^2)^k_*$
\begin{equation}
\label{PsiTheta}
    \widehat\Psi_k(\mu_{\bb})
  = \sum_{n=1}^k \sum_{\Par_n^{k,k}\ni\calT\dashv\bb} \calP_{\calT,\bb} .
\end{equation}
\end{prop}

Proposition \ref{prop:Psi(mu)} is a complicated combinatorial fact based on Proposition \ref{prop:Phi_k}. Note that the equation (\ref{PsiTheta}) looks unexpectedly simple, at least, much more convenient than equation (\ref{Phi=sum+}) for $\Phi_k(\xi)$, $\xi\in 2\pi i\mZ$.

\subsection{Step 6: factorization of the operators $\bW_\calT^b$}

For any $\calT = \{A'_1,A''_1,\ldots,A'_n,A''_n\}\in\Par_n^{k,\kappa}$ we define
$$
    \calT(s)
  = \{A'_s,A''_s\},  \quad
    \alpha_{\calT,\bb}
  = \left\{
    \begin{array}{cc} \prod_{s=1}^n \alpha_{A'_s,A''_s,\bb} &\mbox{if } \bb\vdash\calT, \\
         0                                                  &\mbox{if } \bb\not\vdash\calT.
    \end{array}
    \right.
$$

For any monomial
$$
  \widehat\mu^{\bb} = \widehat\mu^{\beta^{(1)}}\ldots\widehat\mu^{\beta^{(k)}} , \qquad
  \bb\in(\mZ_+^2)^k_*
$$
we define its (multi)degree: $\Deg\widehat\mu^{\bb} = \langle\bb\rangle$. Then
$\Deg \widehat\mu^{\bb'} \widehat\mu^{\bb''}
  = \Deg \widehat\mu^{\bb'} + \Deg \widehat\mu^{\bb''}$ for any $\bb',\bb''\in (\mZ^2_+)^k_*$.
We expand $\bW$ into $\Deg$-homogeneous polynomials:
$$
  \bW = \sum_{b\in\mZ_+^2} \bW^b, \qquad  \Deg \bW^b = b.
$$
For any $b\in\mZ_+^2$, $\calT = \calT(1) \in\Par_1^{k,\kappa}$, and $1\le s\le n$ we put
$$
     H^{\bb_\calT}
  =  \prod_{j\in A'\cup A''} H_{\beta^{(j)}}, \quad
     \bW_\calT^b
  =  \sum_{(\mZ^2_+)^\kappa_*\ni\bb\vdash\calT,\,\langle\bb\rangle_\calT=b}
         \frac{H^{\bb_\calT} \alpha_{\calT,\bb}}{\kappa!}
           \widehat\mu^{\langle\bb\rangle_{\calT}+\one}.
$$

By using Proposition \ref{prop:Psi(mu)} we prove the following statement.

\begin{prop}
For any $b\in\mZ_+^2$
\begin{equation}
\label{Wb=sumsum}
    \bW^b
  = \sum_{1\le n\le k}
    \sum_{\calT\in\Par_n^{k,k}}
    \frac{\kappa_1!\ldots\kappa_n!}{k!}
    \sum_{b_1+\ldots+b_n=b}
      \symm\Big( \bW_{\calT(1)}^{b_1}, \ldots , \bW_{\calT(n)}^{b_n} \Big),
\end{equation}
where $\kappa_s = \#\calT(s)$.
\end{prop}

Any operator $\bW^b_{\{A',A''\}}$, $\{A',A''\}\in\Par_1^{k,\kappa}$, $b\in\mZ_+^2$ depends only on $\# A'$, $\# A''$, and $b$:
\begin{equation}
\label{W=wmu}
  \bW^b_{\{A',A''\}} = w^b(\# A',\# A'') \, \widehat\mu^{b+\one}.
\end{equation}
Therefore by (\ref{Wb=sumsum}) and (\ref{W=wmu})
\begin{eqnarray}
\nonumber
     \bW^b
 &=& \sum_{1\le n\le\kappa\le |b|/2} \bW_{n,\kappa}^b, \\
\label{Wbnkappa}
     \bW^b_{n,\kappa}
 &=& \sum_{\kappa_1+\ldots+\kappa_n=\kappa,\, b_1+\ldots+b_n=b}
          w^{b_1}(\kappa_1) \ldots w^{b_n}(\kappa_n)
          \symm\big( \widehat\mu^{b_1+\one}, \ldots , \widehat\mu^{b_n+\one}
               \big) .
\end{eqnarray}
where
$$
    w^b(\kappa)
  = \sum_{\kappa'+\kappa''=\kappa} \frac{\kappa!}{\kappa'! \kappa''!} w^b(\kappa',\kappa'').
$$

\subsection{Step 7: compute $\sum_\kappa \bW_{n,\kappa}^b$}

For any $s\in\mN$ and $b\in\mZ_+^2$ we put (compare with (\ref{Sigmadk})--(\ref{cond_Sigmadk}))
$$
    \widehat\Sigma_d^{(0,0)}
  = \left\{ \begin{array}{cc}  1 & \mbox{if } d=0, \\
                               0 & \mbox{if } d>0,
            \end{array}
    \right. \quad
    \widehat\Sigma^b_\kappa
  = \frac1{(2\pi)^\kappa}
     \sum_{\langle\bb\rangle=b,\,|\beta^{(j)}|\ge 3} H^\bb, \qquad
    \bb\in(\mZ_+^2)^\kappa.
$$
We put
$$
    p_d
  = \sum_{s=1}^d \frac{(d+s)!}{s!(d+1)!} \widehat\Sigma^{d\cdot\one}_s .
$$

\begin{prop}
\label{lem:sumw}
For any $b\in\mZ_+^2$ and $n,k\in\mN$
$$
     \sum_{\kappa\ge 1} \bW^b_{n,\kappa}
  =  \sum \prod_{j=1}^n
       \frac{2\pi (d_j+1)!}{\sigma_j! (d_j+1-\sigma_j)!}
           \widehat\Sigma_{\sigma_j}^{b_j-d_j\cdot\one} p_{d_j}\,
           \symm\big( \widehat\mu^{b_1+\one},\ldots,\widehat\mu^{b_n+\one} \big),
$$
where the summation is performed under the conditions
$$
  b_1+\ldots+b_n=b , \quad d_1,\ldots,d_n,\sigma_1,\ldots,\sigma_n \ge 1.
$$
\end{prop}

\section{Step 8: equation $\bW=0$ is equivalent to (\ref{eq_conj})}

The equation $\bW=0$ is equivalent to the infinite collection of equations
\begin{equation}
\label{W^b=0}
  \bW^b = 0, \qquad   b\in\mZ_+^2 .
\end{equation}

\begin{prop}
\label{prop_p}
The system (\ref{W^b=0}) is equivalent to the system
\begin{equation}
\label{p_d=0}
  p_d = 0, \qquad   d = 1,2,\ldots
\end{equation}
\end{prop}

Proposition \ref{prop_p} implies Theorem \ref{theo:main}.


\begin{thebibliography}{00}
\bibitem{BM2003} S.~Bolotin and R.~MacKay, Isochronous potentials, in {\it Proc. of the 3rd conference on localization and energy transfer in nonlinear systems}, Luis V\' azquez (ed.) et al., World Sci., 2003, 217--224.
\bibitem{Calogero_book}   F.~Calogero Isochronous systems. Oxford University Press 2008, 264 p.
\bibitem{Z13} G.~Gorni and G.~Zampieri. Global isochronous potentials. Qual. Theory Dyn. Syst. 12 (2013), no. 2, 407–416.
\bibitem{Elf} N.~Elfimov. Master thesis. 2021, Mech. Math Department of Moscow State University.
\bibitem{Still} F.~Stillinger and D.~Stillinger. Pseudoharmonic oscillators and inadequacy of semiclassical quantization. J. Phys. Chem. {\bf 93} (1989) 6890--6892.
\bibitem{Tre_PhD} D.~Treschev, Billiard map and rigid rotation, Phys. D, 255 (2013), 31–34.
\bibitem{Tre_2015} D.V.~Treschev, On a Conjugacy Problem in Billiard Dynamics, Proc. Steklov Inst. Math., 289 (2015), 291–299.
\bibitem{Tre_2017} D.~Treschev, A locally integrable multi-dimensional billiard system, Discrete Contin. Dyn. Syst. Ser. A, 37:10 (2017), 5271–5284.
\bibitem{Tre_2018} V.~Schastnyy, D.~Treschev, On Local Integrability in Billiard Dynamics,
    Exp. Math., 28:3 (2019), 362–368
\end{thebibliography}
\end{document}